\documentclass[12pt,reqno]{amsart}

\usepackage{graphics}
\usepackage{xypic}
\usepackage{amsmath,amsthm}
\usepackage{amscd}
\usepackage{amssymb}
\usepackage{euscript}

\numberwithin{equation}{subsection}

\newcommand{\sqsp}{\renewcommand{\baselinestretch}{1.2}\tiny\normalsize}

\newtheorem{thm}[subsection]{Theorem}

\newtheorem{prop}[subsection]{Proposition}
\newtheorem{cor}[subsection]{Corollary}

\theoremstyle{definition}
\newtheorem{definition}[subsection]{Definition}
\newtheorem{example}[subsection]{Example}
\newtheorem{remark}[subsection]{Remark}

\newcommand{\bK}{\mathbf{k}}

\newcommand{\cP}{\mathcal{P}}

\newcommand{\HP}{H_{\mathcal{P}}}

\newcommand{\CH}{\mathbb{H}}
\newcommand{\ilim}{\stackrel{\mbox{lim}}{\longleftarrow}}
\newcommand{\g}{{\cP^\text{!`}}^{(1)}}

\DeclareMathOperator{\Id}{Id} \DeclareMathOperator{\Hom}{Hom_{\bK}}
\DeclareMathOperator{\End}{End}

\begin{document}

\title{Versal deformation theory of algebras over a quadratic operad}
\author{Alice Fialowski, Goutam Mukherjee and Anita Naolekar}

\date{}

\address{Alice Fialowski, Institute of Mathematics, Eotvos Lorand University, Pazmany
Peter setany 1/C, H-1117 Budapest, Hungary}
\email{fialowsk@cs.elte.hu}

\address{Goutam Mukherjee, Stat-Math Unit, Indian Statistical Institute, 203, B. T. Road, Calcutta- 700108,
India}
\email{goutam@isical.ac.in}

\address{Anita Naolekar, Stat-Math Unit, Indian Statistical Institute, 8th Mile, Mysore Road, RVCE Post, Bangalore- 560059, India}
\email{anita@isibang.ac.in}

\footnote
{The work was supported by INSA and by the Hungarian OTKA grant K77757.}

\begin{abstract} We develop deformation theory of algebras over
quadratic operads where the parameter space is a commutative local
algebra. We also give a construction of a distinguised deformation of an algebra over a quadratic operad with a complete local algebra as its base--the so called `versal deformation'--which induces all
other deformations of the given algebra.
\end{abstract}

\maketitle
\begin{center}
\emph{\large In memory of our good friend Jean-Louis Loday}
\end{center}

\sqsp

\section{Introduction}
Formal one-parameter deformation theory for algebras was originally
introduced for associative algebras by M. Gerstenhaber in the 60'-s
(see \cite{ger, ger2}). Since then it has been applied to many other
algebraic categories. Most of these cases turned out to be algebras
over a suitable quadratic operad. Gerstenhaber's theory was
generalized to an algebra over a binary quadratic operad by D.
Balavoine in 1997 (see \cite{bal2}), and further developed by M.
Kontsevich and Y. Soibelman (see \cite{ks}). In \cite{LV}, J. L. Loday
and B. Vallette has introduced deformation complex for any
$\cP$-algebra $A$ over a quadratic operad $\cP$. It turns out that
this deformation complex is the resulting dg Lie algebra obtained
from convolution dg Lie algebra $\mathfrak{g}_{\cP, A}=
(\mbox{Hom}_{\mathbb S}(\cP^\text{!`}, \mbox{End}_A)), [,])$, where $\cP^\text{!`}$
is the Koszul dual co-operad of $\cP$, twisted by the algebra
structure.

Classical deformation theory is not general enough to describe all
nonequivalent deformations of a given object. To take care of this, one needs to
enlarge the base of deformations from one parameter power series ring to a
local commutative algebra, or more generally, to a complete local
algebra. It is known that under certain cohomological restrictions there exists a "characteristic" versal deformation of a type of algebras, with complete local algebra base, which induces all nonequivalent
deformations and is universal at the infinitesimal level (see e.g.
\cite{sch}). An explicit construction of versal deformation has been given for Lie
\cite{fi1, ff}, associative, infinity \cite{fp} and Leibniz algebras
\cite{fmm}.

The aim of this paper is to show that the above method of construction of versal deformations of a particular type of algebras can be extended to construct versal deformation of an algebra over any quadratic operad. For this it is necessary to develop relevant operadic tools. To this end we develop an obstruction theory for extending a given deformation of an algebra over a quadratic operad with a finite dimensional local commutative algebra base $\mathcal{R}$ to a deformation over  a suitable extension of $\mathcal{R}$ using the operadic calculus as developed in chapter\,12 of \cite{LV}. 

The paper is organized as follows: In Section\,2 we briefly review the basic definitions and results involving algebras over a quadratic operad and the associated deformation complex necessary for our purpose and fix notations which we follow throughtout the text. In Section\,3 we develop
deformation theory of algebras over an operad with local commutative algebra base. In Section\,4 we
study infinitesimal deformations and their properties. Given a deformation $\lambda$ of a $\cP$-algebra $A$, with a finite dimensional base $\mathcal R$ and a local abelian extension $0\rightarrow M\rightarrow \mathcal R'\rightarrow \mathcal R\rightarrow 0$ of $\mathcal R$ by an $\mathcal R$-module $M$, $\mbox{dim}_{\mathcal R}M< \infty$, we address in Section\,5, the question of extending $\lambda$ to a deformation of $A$ with base $\mathcal R'$. It turns out that in general this extension is not always possible. The associated obstruction can be interpreted as a $2$-cochain in the deformation complex associated to $A$. We prove that obstructions are $2$-cocycles, and that vanishing of the associated cohomology class is a necessary and sufficient condition for extending a given deformation to a deformation over a larger base. The results of Section\,5 are crucial in order to construct a versal deformation of an algebra over a quadratic operad which is the content of the next section. In Section\,6 we introduce the notion of formal deformation of algebras over a quadratic operad with a complete local algebra
base and define the notion of versal deformation. Finally we give a
construction of a versal deformation of an algebra over a quadratic operad, using the results developed in Sections\,4 and 5.
\vskip 5em

\section{Preliminaries on operads and operadic cohomology}
\label{sec:prelim}
\bigskip

In this section we recall some basic definitions and results
about (algebraic) operads and the deformation complex of an algebra over a quadratic
operad as in \cite{LV}. General references for operads and related results are \cite{gj,gk,LV,markl,mss,may1,may2}.

The symbol $\mathbb {N}$ denotes the
set of non-negative integers.  Throughout this paper, we work over a
fixed field $\bK$ of characteristic zero. We denote the category of vector spaces over $\bK$ by $\bf{Vect}$. The tensor product of vector
spaces over $\bK$ is denoted by $\otimes$. For any positive integer $n$,
${\mathbb S}_n$ denotes the group of permutations on $n$ elements.
 For any map $f: \otimes_{i=1}^n E_i
\rightarrow F$, $f(x_1,\cdots, x_n)$ will stand for $f(x_1\otimes
\cdots \otimes x_n).$

\subsection{$\mathbb S$-module}
\begin{definition}
An ${\mathbb S}$-module over $\bK$ is a family $$M= \{M(0), M(1), \cdots, M(n),\cdots\}$$ of right $\bK[{\mathbb S}_n]$ modules $M(n)$. An ${\mathbb S}$-module $M$ is finite dimensional if $M(n)$ is finite dimensional for all $n$. A morphism $f:M\longrightarrow N$ between two ${\mathbb S}$-modules $M$ and $N$ is a family of maps $f(n): M(n) \longrightarrow N(n)$ which are ${\mathbb S}_n$-equivariant for all $n$. \end{definition}

Note that to every ${\mathbb S}$-module $M$ is associated a functor, called the $\it {Schur~ functor}$ $\widetilde{M}: \bf {Vect} \longrightarrow \bf {Vect}$ such that $\widetilde{M}(V)= \oplus_{n\geq 0}M(n)\otimes_{{\mathbb S}_n}V^{\otimes n}$, where the left action of ${\mathbb S}_n$ on $V^{\otimes n}$ is given by $\sigma (v_1, v_2, \cdots, v_n)= (v_{\sigma^{-1}(1)}, v_{\sigma^{-1}(2)}, \cdots, v_{\sigma^{-1}(n)})$.

{\it Tensor product} of two $\mathbb S$-modules is the $\mathbb S$-module $M\otimes N$ defined by
$$(M\otimes N)(n)= \oplus_{i+j=n}Ind^{{\mathbb S}_n}_{{\mathbb S}_i\times {\mathbb S}_j} M(i)\otimes N(j).$$

The {\it Composite} of two $\mathbb S$-modules $M$ and $N$ is the $\mathbb S$-module
$$ M\circ N= \oplus_{n\geq 0} M(n)\otimes N^{\otimes n},$$ where $N^{\otimes n}$ is the tensor product of $n$ copies of the $\mathbb S$-module $N$. The category $({\mathbb S}$-Mod, $\circ$, $I$) is a monoidal category with respect to this composite product, where $I$ is the $\mathbb S$-module $(0, \bK, 0, \cdots)$.

Recall that the composite of the two $\mathbb S$-modules $M$ and $N$ satisfies the formula $$\widetilde{(M\circ N)}= \widetilde{M}\circ \widetilde{N},$$ where on the right hand side the symbol $\circ$ stands for the composition of functors.

\subsection{Operads and cohomology}
\begin{definition}
A \emph {symmetric} \emph{operad} $\cP= (\cP, \gamma, \eta)$ is a monoid in the monoidal category $({\mathbb S}$-Mod, $\circ$, $I$). Explicitly, a symmetric operad $\cP= (\cP, \gamma, \eta)$ is an $\mathbb{S}$-module $\cP= \{\cP(n)\}_{n\geq 0}$ endowed with morphisms of $\mathbb{S}$-modules 
$$\gamma: \cP\circ \cP\longrightarrow \cP$$ called composition map and 
$$\eta: I\longrightarrow \cP$$ called the unit map, which makes $\cP$ into a monoid.
\end{definition}

Throughout the paper, by an {\it operad} we shall mean a {\it symmetric operad}.
\begin{definition}
Let $\cP, \mathcal{Q}$ be operads. A {\it morphism of operads} from $\cP$ to $\mathcal{Q}$ is a morphism of $\mathbb{S}$-modules $\alpha: \cP\longrightarrow \mathcal{Q}$, which is compatible with the composition maps.
\end{definition}

\begin{example}
Let $V$ be a vector space over $\bK$ and for every
$n\in {\mathbb N}$ let $\mbox{End}_V(n)= \Hom(V^{\otimes n},
V)$. Then $\mbox{End}_V=\{\mbox{End}_V(n),~n\in \mathbb{N}\}$
is naturally a symmetric operad, called the \emph{endomorphism operad of
$V$}. 
\end{example}

\begin{definition}
Let ${\mathcal P}$ be an operad.  A \emph{$\mathcal P$-algebra}
or an algebra over $\mathcal P$, denoted by $(A,\alpha)$, is a vector space $A$
over $\bK$, equipped with a morphism of operads
$\alpha\colon \mathcal P\rightarrow \End_A$.

A \emph {morphism of $\mathcal P$-algebras} $\phi : (A,\alpha)\longrightarrow
(B,\beta)$ is a $\bK$-linear map $\phi : A \longrightarrow B$ such that for
any $a_1, \cdots,a_n \in A$ and $\mu \in {\mathcal P}(n)$,
$$\phi(\alpha(\mu)(a_1,\cdots, a_n)) = \beta(\mu)(\phi (a_1),\cdots, \phi
(a_n)).$$
\end{definition}


\begin{definition}\label{freeoperad}
The {\it free operad} over the $\mathbb{S}$-module $M$ is an operad $\mathcal{F}(M)$ equipped with an $\mathbb{S}$-module morphism $\eta(M): M\longrightarrow \mathcal{F}(M)$ which satisfies the following universal condition: any $\mathbb{S}$-module morphism $f:M\longrightarrow \mathcal{P}$ where $\cP$ is an operad, extends uniquely into an operad morphism $\widetilde{f}: \mathcal{F}(M)\longrightarrow \cP$, that is $\widetilde{f}\circ \eta(M)= f$.
\end{definition}

The functor $\mathcal{F}$ from the category of $\mathbb{S}$-modules to the category of $\bK$-operads is left adjoint to the forgetful functor from the category of $\bK$-operads  to $\mathbb{S}$-modules. An explicit construction of the free operad can be found in \cite{bjt}.

\begin{definition}
Let $\cP$ be an operad. An ideal of $\cP$ is a sub- $\mathbb{S}$-module $\mathcal{I}$ of $\cP$ such that the operad structure of $\cP$ passes on to the quotient $\cP/\mathcal{I}$.
\end{definition}

Let $E$ be an $\mathbb S$-module and $R$ be a sub-${\mathbb S}$-module such that $R\subseteq {\mathcal F}(E)^{(2)}$, where ${\mathcal F}(E)^{(2)}$ is the graded sub $\mathbb S$-module of the free operad $\mathcal{F}(E)$, which is spanned by the composites of two elements of $E$ (Section\, 5.5.3, \cite{LV}). Such a pair $(E, R)$  is called a quadratic data.

\begin{definition}
Given a quadratic data $(E, R)$, let 
$$\mathcal P(E, R)= \mathcal F(E)/ (R)$$ be the quotient of the free operad $\mathcal F(E)$ over $E$  by the operadic ideal $(R)$ generated by $R$. Then the operad $\mathcal P(E,R)$ is called the quadratic operad associated to the pair $(E, R)$.
\end{definition}

\begin{definition}
An $\mathbb S$-module ${\mathcal C}$ is a {\it cooperad} if it is a comonoid in the monoidal category $(\mathbb S\mbox{-mod}, \circ, I)$. A {\it cofree cooperad} on an $\mathbb S$-module $M$ is the cooperad ${\mathcal F}^c(M)$, which is cofree in the category of conilpotent cooperads.
\end{definition}

Recall that the Koszul dual cooperad of a quadratic operad is defined as follows:
\begin{definition}
The {\it quadratic cooperad} $\mathcal C(E, R)$ associated to the quadratic data $(E, R)$ is the sub-cooperad of the cofree cooperad $\mathcal F^c(E)$ which is universal among the sub-cooperads of $\mathcal F^c(E)$ such that the following composite is $0$, 
$$ \mathcal C \hookrightarrow \mathcal F^c(E)\twoheadrightarrow\mathcal F^c(E)^{(2)}/(R).$$ The {\it Koszul dual cooperad} of the quadratic operad $\cP(E, R)$ is the quadratic cooperad 
$$\mathcal P^\text{!`}:= \mathcal C(sE, s^2 R),$$ where $sE$ denotes the $\mathbb S$-module $E$ whose degree is shifted by $1$.
\end{definition}

We recall the definition of the deformation complex of an algebra over a quadratic operad, \cite{LV}.

Let $A$ be a $\cP$-algebra, where $\cP= \cP(E, R)$ is a quadratic operad. Let  $$\mathfrak{g}=\mathfrak{g}_{\cP^, A}= (\mbox{Hom}_{\mathbb S}(\cP^\text{!`}, \mbox{End}_A), [ , ], \partial)$$ be the {\it convolution dg Lie algebra} associated to $A$. Since $A$ is concentrated in degree $0$, the cohomological degree on the differential graded Lie algebra $$\mathfrak{g}:~~
\mbox{Hom}(A,A)\stackrel{\partial}{\longrightarrow} \mbox{Hom}_\mathbb{S}({\cP^\text{!`}}^{(1)}, \mbox{End}_A)\stackrel{\partial}{\longrightarrow} \mbox{Hom}_\mathbb{S}({\cP^\text{!`}}^{(2)}, \mbox{End}_A) \cdots,$$ is induced by the weight grading. Thus $\mbox{Hom}(A,A)$ is the $0th$ cochain module and $\mbox{Hom}_{\mathbb S}({\cP^\text{!`}}^{(n)}, \mbox{End}_A)$ is the $nth$ cochain module.
Moreover, since $\cP$ is homogeneous quadratic, the coboundary map $\partial$ in the above complex is null. 
We note that the set of $\cP$-algebra structures
on a space $A$ is in one-to-one correspondence with the set of
Maurer-Cartan elements of degree $1$ of $\mathfrak{g}$, i.e. all elements $\phi$ in $\mbox{Hom}_\mathbb{S}({\cP^\text{!`}}^{(1)}, \mbox{End}_A)$ satisfying $[\mathbb{\phi},\mathbb{\phi}]=0$.
Given such an element $\phi$ one can define a differential
$\partial_{\phi}$ on the Lie algebra $\mathfrak{g}$, which makes it into a differential graded Lie algebra. The differential $\partial_\phi$ is called the twisted differential and $\mathfrak{g}^\phi=\mathfrak{g}_{\cP, A}^\phi=(\mbox{Hom}_\mathbb{S}(\cP^\text{!`}, \mbox{End}_A)), [,], \partial_{\phi})$ the twisted differential graded Lie algebra. The twisted
differential $\partial_\phi$ is given by: $\partial_\phi(f)= [\phi, f]$. The underlying cochain
complex of this twisted dg Lie algebra is the deformation complex
that we intend to work with. We end this section with the following definition.
\begin{definition} For any $\cP$-algebra $(A, \pi)$ where $\cP$ is a quadratic operad, we define 
$$H^*_{\cP}(A): = H_*(\mathfrak{g}^\pi)= H_*(\mathfrak{g}, \partial_\pi).$$
\end{definition}

\section{Deformations}
\bigskip

Let $\mathcal R$ be a commutative local unital algebra with unit $1_{\mathcal R}$ over $\bK$. Let $\epsilon:
\mathcal R\rightarrow \bK, \epsilon(1_{\mathcal R})=1$ be the canonical
augmentation map, and ${\mathfrak{M}}= ker(\epsilon)$ be the unique
maximal ideal in $\mathcal R$. In this section, we study the notion of
deformation of an algebra over a quadratic operad with base $\mathcal R$ and
 its properties.

Let $\cP=\cP(E, R)$ be a quadratic operad. We will denote by $\cP_\mathcal R$ the operad which is obtained by
extension of $\cP$ to the category of modules over $\mathcal R$, in other words,
$\cP_\mathcal R(n) = \mathcal R\otimes \cP(n)$ for all $n\in \mathbb {N}.$
Let $(A, \pi)$ be a $\cP$-algebra. Let $A_\mathcal R= \mathcal R\otimes A$ denote the
extension of $A$. Then $A_\mathcal R$ can be viewed as a  $\cP_\mathcal R$-algebra  by
extending $\pi \colon \cP \longrightarrow \mbox{End}_A$ to
$\pi_\mathcal R \colon \cP_\mathcal R \longrightarrow \mbox{End}_{A_\mathcal R}$, since
$$\mbox{Hom}_\mathcal R(A_\mathcal R^{\otimes_\mathcal Rn},A_\mathcal R) \cong \mathcal R\otimes \Hom(A^{\otimes
n},A).$$ Moreover, $A=\bK \otimes A$ can be viewed as a $\cP_\mathcal R$-algebra
by considering $A$ as a module over $\mathcal R$ via $\epsilon$, that is $\mathcal R$-module structure  on $A$ is given by $r\cdot a =
\epsilon (r)a$, for $r\in \mathcal R$ and $a\in A$.

\begin{definition}
A deformation $\lambda$ of a $\cP$-algebra $(A, \pi)$ with base $(\mathcal R, \mathfrak{M})$
is a morphism of operads $\lambda: \cP_\mathcal R
\longrightarrow \End_{A_\mathcal R}$ such that  $(\epsilon \otimes \Id):
 A_\mathcal R \longrightarrow \bK \otimes A\cong A$ is a
$\cP_\mathcal R$-algebra morphism.
\end{definition}

In other words, a deformation of $(A, \pi)$ with base $(\mathcal R, \mathfrak{M})$ is a $\cP_\mathcal{R}$-algebra structure on $A_\mathcal{R}$ which reduces to $\pi$, modulo $\mathfrak{M}$: $A\cong \bK\otimes_\mathcal{R}(\mathcal{R}\otimes A)$. 

We have the following result (Proposition\,12.2.6, \cite{LV}).

\begin{prop}\label{MC}
For any $\cP$-algebra $(A, \pi)$, the set of all deformations of $(A, \pi)$ with base $(\mathcal{R}, \mathfrak{M})$ is in one-to-one correspondence with the set of all Maurer-Cartan elements in the dg Lie algebra $\mathfrak{M}\otimes \mathfrak{g}^\pi$.
\end{prop}

\begin{definition}
Suppose $\lambda_1$ and $\lambda_2$ are two deformations of a
$\cP$-algebra $(A, \pi)$ with the same base $(\mathcal R, \mathfrak{M})$. They
are said to be equivalent if there exists a $\cP_\mathcal R$-algebra
isomorphism $\phi: (A_\mathcal R,\lambda_1)\longrightarrow
(A_\mathcal R,  \lambda_2)$ such that
$$\begin{array}{rcl}
A_\mathcal R &\stackrel{\phi}{\longrightarrow} & A_\mathcal R\\
\epsilon\otimes \Id \searrow & & \swarrow \epsilon \otimes \Id\\
&\bK \otimes A&
\end{array}
$$
commutes.
\end{definition}

\begin{definition}
Let $\mathcal R$ be a complete local algebra, $\mathcal R=
\stackrel{\leftarrow}{\mbox{lim}}(\mathcal R/\mathfrak{M}^n)$, $\mathfrak{M}$ denoting the maximal ideal in $\mathcal R$. A formal
deformation of a $\cP$-algebra $(A, \pi)$ with base $(\mathcal R, \mathfrak{M})$
is a $\cP_{\mathcal R}$-algebra structure on the completed tensor product
$\mathcal R \widehat{\otimes}A=
\stackrel{\leftarrow}{\mbox{lim}}((\mathcal R/\mathfrak{M}^n)\otimes A)$, such that $\epsilon \widehat{\otimes}\Id~: \mathcal R
\widehat{\otimes} A\longrightarrow \bK\otimes A \cong A$ is a $\cP_{\mathcal R}$-algebra morphism.
\end{definition}

\begin{example}
Let $\mathcal R= \bK[[t]]$ be the ring of formal power series with
coefficients in $\bK$. Then a formal deformation of a
$\cP$-algebra $(A, \pi)$ over $\mathcal R$ is precisely the formal
`$1$-parameter' deformation of $(A, \pi)$.
\end{example}

\begin{definition}
Let $\lambda$ be a deformation of the $\cP$-algebra $(A, \pi)$ with base
$(\mathcal R, \mathfrak{M})$ and augmentation $\epsilon: \mathcal R\longrightarrow \bK$. Let $\mathcal R'$ be
another commutative local algebra with unit, and augmentation
$\epsilon ': \mathcal R'\longrightarrow \bK$ with $\mbox{Ker}(\epsilon ')= {\mathfrak{M}}'$.
Let $f: \mathcal R\longrightarrow \mathcal R'$ be an
algebra homomorphism with $f(1_{\mathcal R})=1_{\mathcal R'}$. Then  $\epsilon ' \circ f= \epsilon$.

Consider $\mathcal R'$ as an $\mathcal R$-module via the map $f$: $r'\cdot r= r' f(r)$
so that $$\mathcal R'\otimes A= (\mathcal R'\otimes_\mathcal R \mathcal R)\otimes A= \mathcal R'\otimes_\mathcal R (\mathcal R \otimes A).$$
Then the {\it push-out of the deformation} $\lambda$ is the deformation
$f_*\lambda $ of $(A, \pi)$ with base $(\mathcal R', \mathfrak{M}')$, defined by
$$\begin{array}{ll}
&f_*\lambda(1_{\mathcal R'}\otimes \mu)\{r_1' \otimes_\mathcal R (r_1 \otimes a_1), r_2'
\otimes_\mathcal R (r_2 \otimes a_2), \cdots, r_n' \otimes_\mathcal R (r_n \otimes a_n)\}\\
=& r_1'r_2'\cdots r_n'\otimes_\mathcal R
\lambda(1_{\mathcal R}\otimes \mu)(r_1\otimes a_1, r_2\otimes a_2,\cdots, r_n\otimes a_n),
~\mu \in \cP^{\text{!`}(1)}(n).
\end{array}
$$
\end{definition}

It is straightforward to see that $f_* \lambda$ is a deformation of
$(A,\pi)$ with base $(\mathcal R', \mathfrak{M}')$.

\begin{remark}\label{decom}
Note that by the proof of Proposition \ref{MC}, as in \cite{LV},  if $\lambda$ is a deformation of the $\cP$-algebra $(A, \pi)$ with a base $(\mathcal R, \mathfrak{M})$, then $\lambda$ can be expressed as 
\begin{equation}\label{lambda}
\lambda= \pi +\sum_{i=1}^s m_i\otimes \phi_i,
\end{equation}
 where $\pi$ is a Maurer-Cartan element of $\mathfrak{g}^\pi$ and $\sum_{i=1}^s m_i\otimes \phi_i$ is a Maurer-Cartan element of $\mathfrak{M}\otimes \mathfrak{g}^\pi$. Also the pushout $f_* \lambda$ is given by
\begin{equation}
\label{pushout}
f_*\lambda = \pi + \sum_{i=1}^s f(m_i) \otimes \phi_i,
\end{equation} where $\sum_{i=1}^s f(m_i) \otimes \phi_i$ is a Maurer-Cartan element of $\mathfrak{M}'\otimes \mathfrak{g}^\pi$.
\end{remark}

\begin{definition}
A deformation $\lambda$ of $(A, \pi)$ with base $(\mathcal R, \mathfrak{M})$ is called
{\it infinitesimal} if, in addition, ${\mathfrak{M}}^2=0$.
\end{definition}

To consider the equivalence of infinitesimal deformations the
cohomology comes into play naturally.

Let $\lambda$ be an infinitesimal deformation of $(A, \pi)$, with base $(\mathcal R, \mathfrak{M})$ such that $$\lambda= \pi +\sum_{i=1}^s m_i\otimes \phi_i.$$ 
Let $\xi \in {\mathfrak{M}'}= \Hom(\mathfrak{M}, \bK)$. Clearly, $\xi$ can be viewed
as an element of  $\Hom(\mathcal R, \bK)$ with $\xi(1_\mathcal R)=0$. 

Define a $1$-cochain $\alpha_{\lambda, \xi}$ of $\mathfrak{g}^\pi$ by
$$\alpha_{\lambda, \xi}= (\xi\otimes \mbox{Id}) \lambda.$$

\begin{thm}
For any infinitesimal deformation $\lambda$ of a $\cP$-algebra
$(A,\pi)$, $\alpha_{\lambda, \xi}$ is a $1$-cocycle.
\end{thm}
{\bf Proof.} As in remark \ref{decom}, we write $\lambda= \pi +\sum_{i=1}^s m_i\otimes \phi_i.$ 
Then the $1$-cochain $\alpha_{\lambda, \xi}$ can be expressed as
$$\begin{array}{rl}
\alpha_{\lambda, \xi}&= (\xi\otimes \mbox{Id}) \lambda\\
&=(\xi\otimes \mbox{Id})( \pi +\sum_{i=1}^s m_i\otimes \phi_i)\\
&= \sum_{i=1}^s \xi(m_i) \phi_i, 
\end{array}$$ as $\xi(1_{\mathcal R})=0$.

The fact that $\sum_{i=1}^s m_i\otimes \phi_i$ ia a Maurer-Cartan element in $\mathfrak{M}\otimes \mathfrak{g}^\pi$, and the fact that $\mathfrak{M}^2=0$ implies that $$\sum_{i=1}^s m_i\otimes [\phi_i, \pi]=0.$$
Therefore, $$\sum_{i=1}^s\xi(m_i) [\pi, \phi_i]=(\xi\otimes \mbox{Id})(\sum_{i=1}^s m_i\otimes [\pi, \phi_i])=0.$$ 
Using the above observations, we deduce $$\partial_\pi(\alpha_{\lambda,\xi})= [\pi, \alpha_{\lambda, \xi}]= [\pi, \sum_{i=1}^s\xi(m_i)\phi_i]= \sum_{i=1}^s\xi(m_i)[\pi, \phi_i]=0.$$ \qed

Let us define for $\xi \in {\mathfrak{M}}'$ the cohomology class of the
cocycle $\alpha_{\lambda, \xi}$ by $a_{\lambda, \xi}.$
The correspondence $\xi \longmapsto a_{\lambda, \xi}$ defines
a map $$a_\lambda: {\mathfrak{M}}' \longrightarrow \HP^1(A).$$

\begin{thm}\label{equi}
Let $\lambda_1$ and $\lambda_2$ be two infinitesimal deformations of
$(A, \pi)$ with base $(\mathcal R, \mathfrak{M}) $. Assume that $\mathcal{R}$ is of finite dimension.
Then the deformations $\lambda_1$ and $\lambda_2$ are equivalent iff
$\alpha_{\lambda_1, \xi}$ and $\alpha_{\lambda_2, \xi}$ represent the same
cohomology class, that is $a_{\lambda_1, \xi} = a_{\lambda_2, \xi}$ for
$\xi \in \mathfrak{M}'$.
\end{thm}

{\bf Proof.} This result is proved in Theorem 12.2.7, \cite{LV} for $R= \bK[t]/(t^2)$. The present theorem is proved 
by following the same idea. We only give the essential steps. 

Let $\{m_i\}_{1\leq i\leq r}$ be a basis of $\mathfrak{M}$ and $\{\xi_i\}_{1\leq i\leq r}$ be the dual basis of $\mathfrak{M}'$. By definition
$\lambda_1$ and $\lambda_2$ are  equivalent if and only if there exists a
$\cP_\mathcal R$-algebra isomorphism
\begin{equation}
\label{iso}
\rho: A_\mathcal R\longrightarrow A_\mathcal R,~ \text{such that}~
(\epsilon \otimes \Id)\circ \rho= \epsilon \otimes \Id.
\end{equation}

Since $A_\mathcal R= \mathcal R\otimes A= (\bK \oplus {\mathfrak{M}})\otimes A \cong
A \oplus ({\mathfrak{M}}\otimes A)$, the isomorphism  $\rho$ can
be written as $\rho = \rho_1 + \rho_2$ where $\rho_1: A\longrightarrow
A$ and $\rho_2: A \longrightarrow {\mathfrak{M}}\otimes A$.

By compatibility (\ref{iso}), we get $\rho_1=\Id$.
Using the adjunction property of tensor products, we have
$$
\Hom(A;{\mathfrak{M}}\otimes A) \cong {\mathfrak{M}}\otimes \Hom(A,A)
\cong \Hom({\mathfrak{M}}'; \Hom(A,A)),
$$
where the isomorphisms are given by
 \begin{equation}
\label{adiso}
\rho_2\longmapsto \sum_1^r m_i \otimes \phi_i\longmapsto \sum_i^r \chi_i.
\end{equation}
Here $\phi_i = (\xi_i \otimes \mbox{id})\circ \rho_2$, $\chi_i
(\xi_j)= \delta_{i, j}\phi_i$.

Thus we can write,
$$\rho=\mbox{Id}+ \sum_1^r m_i \otimes \phi_i.$$

Recall that, $\rho$ is a $\cP_\mathcal R$-algebra morphism iff
\begin{equation}
\label{rho1}
\rho( \lambda_1 (1_{\mathcal R}\otimes \mu))= \lambda_2(1_{\mathcal R}\otimes \mu)(\rho, \rho, \cdots, \rho),\end{equation}
$\mu \in {\cP^\text{!`}}^{(1)} (n)$, where on the right hand side we have $n$ many copies of $\rho$.

Let us set $\psi_i^k= \alpha_{\lambda_k, \xi_i}$, $i=1,2,\ldots, r$ and
$k=1,2.$ Then we have

\begin{equation}
\lambda_k(1_{\mathcal R}\otimes \mu)= \pi (\mu)+ \sum_1^rm_i \otimes \psi_i^k(\mu)
\end{equation} for all $\mu \in {\cP^\text{!`}}^{(1)} $ .
We explicitly write both sides of the equation (\ref{rho1}), using the fact that the deformations involved are infinitesimal.
$$
\begin{array}{ll}
& \lambda_2 (1_{\mathcal R}\otimes \mu)(\rho,\rho, \cdots, \rho)\\
=&  \pi(\mu) + \sum_1^r m_i \otimes \psi_i^2(\mu) +\sum_{k=1}^n \sum_1^r m_i \otimes( \pi(\mu)( \mbox{Id}, \mbox{Id}, \cdots, \phi_i, \mbox{Id}, \cdots)),
\end{array}
$$ where in the last summand $\phi_i$ is in the $k$th slot.

On the other hand,
\begin{align*}
&\rho(\lambda_1(1_{\mathcal R}\otimes \mu))\\
=& \rho( \pi (\mu)+ \sum_1^r m_i \otimes \psi_i^1(\mu))\\
=& \pi (\mu)
+ \sum_1^r m_i \otimes \phi_i\circ \pi(\mu)+\sum_1^r m_i \otimes \psi_i^1(\mu) 
\end{align*}

It follows from (\ref{rho1}), $$\sum_1^r m_i \otimes (\psi_i^2
- \psi_i^1)+ \sum_1^r m_i \otimes \partial_{\pi} \phi_i=0.$$
Hence,
\begin{align*}
\partial_{\pi}\phi_i=& \psi_i^2- \psi_i^1\\
=&\alpha_{\lambda_2, \xi_i}- \alpha_{\lambda_1, \xi_i}
~~~~~~~~~~~\mbox{for all}~~i=1, \ldots, r.
\end{align*}
So, $a_{\lambda_1, \xi}= a_{\lambda_2, \xi}$, for all $\xi \in {\mathfrak{M}}'$.\qed


Let $\mathcal R$ be a local algebra with maximal ideal $\mathfrak{M} $.
Then $\mathcal R/{\mathfrak{M}}^2$ is local with maximal ideal ${\mathfrak{M}}/{\mathfrak{M}}^2$ and 
${({\mathfrak{M}}/{\mathfrak{M}}^2)}^2=0$. Let
$q\colon \mathcal R\longrightarrow \mathcal R/{\mathfrak{M}}^2$ be the projection map. Let $\lambda$ be a deformation of $(A, \pi)$ with base $(\mathcal R, \mathfrak{M})$. Then the deformation $q_*\lambda$ is infinitesimal and we have a map $$a_{q_* \lambda}: ({\mathfrak{M}}/{\mathfrak{M}}^2)' \longrightarrow \HP^1(A).$$

 \begin{definition}
The dual space $({\mathfrak{M}}/{\mathfrak{M}}^2)'$ is called the {\it tangent space} of $\mathcal R$ and is denoted by $T{\mathcal R}$. The mapping
$$a_{q_* \lambda}: ({\mathfrak{M}}/{\mathfrak{M}}^2)' \longrightarrow \HP^1(A)$$
is called the ${\it differential}$ of $\lambda$ and is denoted by
$d\lambda$. In particular, if ${\mathfrak{M}}^2=0$, then the differential
$d\lambda $ of the infinitesimal deformation $\lambda $ is the map $
a_{\lambda}.$
\end{definition}

\begin{cor}
If two deformations $\lambda_1$ and $\lambda_2$ of a $\cP_\mathcal R$-algebra
$(A, \pi)$ are equivalent, then their ${\it differentials}$ are equal.
\end{cor}
\medskip

\section{Universal Infinitesimal deformation}\label{universal}
\bigskip

In this section we construct a specific example of an infinitesimal
deformation of a $\cP$-algebra satisfying finite dimensionality of the
first cohomology module. We shall also prove a fundamental property of
this deformation. In the last section we will see that this
infinitesimal deformation is the first step of an inductive construction of a versal deformation.

Let $(A,\pi)$ be a given $\cP$-algebra satisfying the condition
$\mbox{dim}~H^1_{\cP}(A) < \infty.$ Let us denote $H^1_{\cP}(A)$ by
$\mathbb {H}.$ Consider the $\bK$-algebra $C_1 = \bK \oplus \mathbb
{H'}$ with the following structure:
$$(k_1,h_1) \cdot (k_2,h_2) = (k_1k_2, k_1h_2+k_2h_1).$$
Clearly, $C_1$ is local with maximal ideal ${\mathfrak{M}}= \mathbb {H'}$ and
${\mathfrak{M}}^2 = 0.$

Fix a homomorphism
$$\sigma \colon \mathbb {H} \longrightarrow \mathfrak{g}_{\cP, A}^1 = \mbox{Hom}_{\mathbb S}({\cP^{\text{!`}}}^{(1)}, \mbox{End}_A)$$
which takes a cohomology class into a representative cocycle. We note
that
$$\mathbb {H'} \otimes A = \Hom(H^1_{\cP}(A), \bK)\otimes A \cong
\Hom(H^1_{\cP}(A),A)$$ and
$$A_{{\mathcal C}_1}= C_1\otimes A \cong ( \bK \oplus \mathbb {H'}) \otimes A \cong A\oplus \Hom
(\mathbb {H}, A).$$

Define a $\cP_{C_1}$-algebra structure $\eta_1$ on $A_{C_1}$ by

\begin{equation}
\label{eta}
\eta_1=  \pi+\sum_{i=1} ^n{h_i}^* \otimes \sigma(h_i)
\end{equation}
as a map $C_1\otimes {\cP^\text{!`}}^{(1)}\longrightarrow (\bK \oplus \mathbb{H}')\otimes \mbox{End}_A= C_1\otimes \mbox{End}_A$, where $\{h_i\}_{1\leq i\leq n}$ denotes a finite basis of $\mathbb{H}$ and $\{{h_i}^*\}_{1\leq i\leq n}$ denotes the dual basis.

\begin{prop}
For any homomorphism $\sigma \colon \mathbb {H} \longrightarrow \mathfrak{g}_{\cP, A}^1 = \mbox{Hom}_{\mathbb S}({\cP^\text{!`}}^{(1)}, \mbox{End}_A)$, $(A_{C_1}, \eta_1)$ is a $\cP_{C_1}$-algebra.
\end{prop}
{\bf Proof.} We need to show that $\eta_1$ is a Maurer-Cartan element in $C_1\otimes \mathfrak{g}^\pi$. The coboundary of the complex $ C_1\otimes \mathfrak{g}^\pi$ is given by 
$$\bar{\partial_\pi}:=\mbox{Id}_{C_1}\otimes \partial_\pi.$$  From $(\ref{eta})$, 
$$ \eta_1=  \pi+\sum_i {h_i}^* \otimes \sigma(h_i).$$
So $$\bar{\partial_\pi}(\eta_1)=\mbox{Id}_{C_1}\otimes [\pi, \eta_1]= \mbox{Id}_{C_1} \otimes [\pi, \pi] + \sum_i h_i^* \otimes [\pi, \sigma(h_i)]\big) =0$$
using the fact that $\sigma(h_i)$ is a $1$-cocycle in $\mathfrak{g}^\pi$.  \qed

\begin{prop}
Up to isomorphism, the $\cP_{C_1}$-algebra structures of $A_{C_1}$ does not
depend on the choice of $\sigma$.
\end{prop}
{\bf Proof.} Let $\sigma': \CH \longrightarrow \mathfrak{g}_{\cP, A}^1$ be another choice of $\sigma$ and denote the corresponding
$\cP_{C_1}$-algebra structure on $C_1\otimes A$ by $\eta'$. Then for $h \in
\CH$, $\sigma(h)$ and $\sigma'(h)$ are two $1$-cocycles
of $A$, representing the same cohomology class, that is, 
$\sigma(h) - \sigma'(h)$ is a $1$-coboundary. Let
$\sigma(h) - \sigma'(h)= \partial_\pi(\gamma(h))$ where
$\gamma: \CH\longrightarrow \Hom (A, A)$. Here $\Hom (A, A)$ is the $0$ cochains of $\mathfrak{g}^\pi.$ 
Using the identification $$C_1\otimes A\cong A \oplus \Hom(\CH, A)$$
define a $C_1$-linear automorphism
$\rho: C_1\otimes A \longrightarrow C_1\otimes A$ by $\rho(a,\phi)= (a, \bar{\phi})$ where $\bar{\phi}(h)= \phi(h) + \gamma
(h)(a)$. Need to show that $\rho$ preserves the $\cP_{C_1}$-algebra structure,
that is,
\begin{equation}
\label{rho}
\rho(\eta _1 (1_{C_1}\otimes\mu))
=\eta'(1_{C_1}\otimes\mu)\circ(\rho, \rho, \cdots, \rho),
\end{equation} where the number of copies of $\rho$ on the right hand side is same as the (homogeneous) degree of $\mu$ in $\g$.
Now observe that using the isomorphism 
$$\Hom(\mathbb{H}, A)\cong \Hom(\mathbb{H}, \bK)\otimes A,$$ we can rewrite the expression (\ref{eta}) of $\eta_1$  as
$$\eta_1(1_{C_1}\otimes \mu)((a_1, \phi_1), \cdots, (a_n, \phi_n))= (\pi(\mu)(a_1, \cdots, a_n), \psi_\mu)$$ where
$\psi_\mu(h)= \sigma(h)(\mu)(a_1, \cdots, a_n) + \sum_i \pi(\mu)(a_1, \cdots, \phi_i(h),\cdots, a_n),$ where in the summation, $\phi_i$ is at the $i$th slot.
Similarly, $$\eta'(1_{C_1}\otimes \mu)((a_1, \phi_1), \cdots, (a_n, \phi_n))= (\pi(\mu)(a_1, \cdots, a_n), \psi'_\mu)$$
 where
$$\begin{array}{rl}\psi'_\mu(h)=& \sigma'(h)(\mu)(a_1, \cdots, a_n) + \sum_i \pi(\mu)(a_1, \cdots, \phi_i(h),\cdots, a_n)\\
=&(\sigma(h)- \partial_\pi(\gamma(h))) (\mu)(a_1, \cdots, a_n) + \sum_i \pi(\mu)(a_1, \cdots, \phi_i(h),\cdots, a_n).
\end{array}$$

We evaluate both sides of (\ref{rho}) using the expressions of $\eta_1$, $\eta'$ as given above and the definitions of $\rho$ and the coboundary $\partial_\pi$, to get $\psi'_\mu(h)- \psi_\mu(h)=0$, for all $h\in \mathbb{H}$. Then equation (\ref{rho}) holds.\qed

The main property of the infinitesimal deformation $\eta_1$ is its
universality in the class of infinitesimal deformations with finite
dimensional base.

\begin{thm}\label{couniversal}
For any infinitesimal deformation $\lambda$ of the $\cP$-algebra $A$ with a
finite dimensional local base $(\mathcal R, {\mathfrak{M}})$, there exists a unique homomorphism
$\phi: C_1\rightarrow \mathcal R$ such that $\lambda$ is equivalent to the push-out
$\phi_*\eta_1$.
\end{thm}
{\bf Proof.} Let $a_{\lambda}: {\mathfrak{M}}' \rightarrow \HP^1(A)=\CH$ denote
the differential of $\lambda$,
$$a_{\lambda}: \xi \mapsto a_{\lambda, \xi}=[\alpha_{\lambda,
\xi}],~~\xi \in {\mathfrak{M}'}.$$

Consider the map
$\phi= \Id \oplus a_{\lambda}': \bK \oplus \CH'\rightarrow \bK \oplus \mathfrak{M}
= \mathcal R$. Let $\{m_i\}_{i=1}^r$ be a basis of $\mathfrak{M}$ and
$\{\xi_i\}$ be the dual basis.

It is enough to show $\alpha_{\phi_* \eta_1}= \sigma\circ a_\lambda$
(Theorem \ref{equi}) . Let $\{h_1, \cdots, h_n\}$ be a basis of $\CH$, and
$\{{h_1}^*, \cdots, {h_n}^*\}$ be the corresponding dual basis of $\CH '$.

Then by Remark \ref{decom}, equation\,(\ref{pushout}) we have
$$
 \phi_*\eta_1=  \pi+ \sum_{j=1}^n \phi({h_j}^*) \otimes \sigma(h_j).$$

Now,
\begin{align*}
a_\lambda '({{h_j}^*})&=\sum_{i=1}^r \xi_i(a_\lambda '({h_j}^*))m_i, \\
\mbox{and} \hspace{40pt} a_\lambda(\xi_i)&=\sum_{j=1}^n {h_j}^*(a_\lambda(\xi_i))h_j.
\end{align*}
Thus
\begin{align*}
&\alpha_{\phi_*\eta_1}(\xi_i)(\mu; a_1, \cdots, a_n)\\
=&(\xi_i \otimes \Id) \phi_*\eta_1(\mu)(1_\mathcal R\otimes a_1,\cdots, 1_\mathcal R\otimes a_n)\\
=&(\xi_i \otimes \Id)\Big(\pi (\mu)(a_1,\cdots, a_n)
+ \sum_{j=1}^n \phi({h_j}^*)\otimes \sigma(h_j)(\mu; a_1, \cdots, a_n)\Big)\\
=&(\xi_i \otimes \Id)(\sum_{j=1}^n a_\lambda '({h_j}^*)\otimes \sigma(h_j)(\mu; a_1,\cdots, a_n))\\
=&\sum_{j=1}^n \xi_i(a_\lambda '({h_j}^*))\otimes \sigma(h_j)(\mu; a_1, \cdots, a_n)\\
=&\sum_{j=1}^k {h_j}^*(a_\lambda (\xi_i))\otimes \sigma(h_j)(\mu; a_1, \cdots, a_n)\\
=&\sigma\big(\sum_{j=1}^n {h_j}^*(a_\lambda(\xi_i))h_j\big)(\mu; a_1, \cdots, a_n)\\
=&\sigma \circ a_\lambda(\xi_i)(\mu; a_1,\cdots,  a_n)
\end{align*}

Therefore, $\alpha_{\phi_*\eta_1}= \sigma\circ a_\lambda$.
The uniqueness of the homomorphism $\phi: C_1\longrightarrow \mathcal{R}$ follows from the definition of $\phi$ and the fact that two infinitesimal deformations $\theta$ and $\theta'$ are equivalent if and only if the corresponding maps $a_{\theta}$ and $a_{\theta'}$ are equal, (cf. Theorem\,\ref{equi}).\qed

\medskip
\section{Deformation Extensions}
\bigskip

Let us recall some definitions and results from \cite{harrison}. Let $\mathcal R$ be a commutative algebra over $\bK$ . Let $(C_q(\mathcal R), \delta)$ denote the standard Hochschild complex, where $C_q(\mathcal R)$ is the $\mathcal R$-module $\mathcal R^{ \otimes (q+1)}$ with $\mathcal R$ acting on the first factor by multiplication of $\mathcal R$. Let $Sh_q(\mathcal R)$ be the $\mathcal R$-submodule of $C_q(\mathcal R)$ generated by chains
\begin{equation*}
\begin{split}
&s_p(r_0,r_1,r_2,\cdots,r_q)\\
=&~\sum_{\sigma\in Sh(p,q-p)}    {sgn(\sigma)(r_0,r_{\sigma^{-1}(1)},r_{\sigma^{-1}(2)},\cdots,r_{\sigma^{-1}(q)}) \in C_q (\mathcal R)}
\end{split}
\end{equation*}
$\mbox{for}~ r_1,r_2,\cdots,r_q \in \mathcal R ~;~0<p<q $.

Then $Sh_{*}$ is a subcomplex of $C_{*}(\mathcal R)$ and hence we have a complex called the {\it Harrison complex}
\[Ch_*(\mathcal R)=\{ Ch_q(\mathcal R),\delta\}~;~ Ch_q(\mathcal R)=C_q(\mathcal R)/Sh_q(\mathcal R).\]
For an $\mathcal R$-module $M$, the Harrison cochain complex defining the Harrison cohomology with coefficients in $M$ is given by
$Ch^{*}(\mathcal R~;M)=Hom_{\mathcal R}(Ch_{*}(\mathcal R), M)$.
\begin{definition}
 For an $\mathcal R$-module $M$ we define
$$ H_{Harr}^{q}(\mathcal R;~M)= Z^q(\mathcal R, M)/ B^q(\mathcal R, M)= H^q(Ch^*(\mathcal R, M)), $$ where $Z^q(\mathcal R, M)$ and $B^q(\mathcal R, M)$ are the spaces of $q$-cocycles and $q$-coboundaries respectively.
\end{definition}

\begin{prop}\label{coefficients in M}
Let $\mathcal R$ be a commutative local algebra with maximal ideal $\mathfrak{M}$, and let $M$ be an $\mathcal R$-module with $\mathfrak{M}M=0$.  Then we have the canonical isomorphism
$$ H_{Harr}^{q}(\mathcal R;~M)\cong  H_{Harr}^{q}(\mathcal R; \bK)\otimes M.$$
\end{prop}

\begin{definition}
A {\it (split) abelian extension} or {\it square-zero extension} $\mathcal R'$ of $\mathcal R$ by an $\mathcal R$-module $M$ is a $\bf k$-algebra $\mathcal R'$ together with an exact sequence of $\bf k$-modules
$$0\rightarrow M\stackrel{i}{\rightarrow} \mathcal R'\stackrel{p}{\rightarrow}\mathcal R\rightarrow 0,$$
where $p$ is an algebra homomorphism so that $N=i(M)$ is an $\mathcal R'$-module and this  $\mathcal R'$-module structure is induced by the $\mathcal R$-module structure on $M$ as $r'i(m) = i(p(r')m).$ In particular, $N$ is an ideal in $\mathcal R'$ satisfying $N^2 =0$.
\end{definition}

\begin{definition} An abelian extension $\mathcal R'$ of $\mathcal R$ by $\mathcal R$-module $M$ is called a {\it local abelian extension} if in addition, $\mathfrak{M}M=0$, where $\mathfrak{M}$ is the maximal ideal of $\mathcal R$.
\end{definition}

{\it Henceforth, by an extension we shall always mean a local abelian extension.}

\begin{remark}
Note that as $\mathcal R$ is local, $\mathcal R'$ is also local with ${\mathfrak{M}}_{\mathcal R'}=p^{-1}(\mathfrak{M})$ as its maximal ideal. Moreover, the condition ${\mathfrak{M}}M=0$ clearly implies that for any $x\in \mathfrak{M}_{\mathcal R'}$ and $n\in N,~~ xn =0$.
\end{remark}
We will use the following results relating
Harrison cohomology and extensions of the algebra $\mathcal R$ by means of $M$, \cite{harrison}.

\begin{prop}\label{cohomology class corresponds to extension}
\begin{enumerate}
\item[(i)] The space $ H_{Harr} ^1 (\mathcal R ;M)$ is isomorphic to the space of derivations $\mathcal R \longrightarrow M$.
\item[(ii)] Elements of $ H_{Harr} ^2 (\mathcal R;M)$  correspond bijectively to isomorphism classes of extensions
$$ 0\longrightarrow M\longrightarrow \mathcal R' \longrightarrow \mathcal R\longrightarrow 0 $$
of the algebra $\mathcal R$ by means of $M$.
\item[(iii)] The space $H_{Harr}^1(\mathcal R;M)$ can also be interpreted as the group of automorphisms of any given  extension
 of $\mathcal R$ by $M$.
\end{enumerate}
\end{prop}

\begin{cor}\label{TA}
 If $\mathcal R$ is a local algebra with the maximal ideal $\mathfrak{M}$, then $$H_{Harr} ^1(\mathcal R~;\bK)\cong \Big(\frac{\mathfrak{M}}{\mathfrak{M}^2}\Big)'=T\mathcal R.$$
\end{cor}

Let $\mathcal R$ be a finite dimensional commutative, unital, local algebra with
augmentation $\epsilon$, and maximal ideal $\mathfrak{M}$. Let $\lambda$ be a deformation of a $\cP =\cP(\bK,E, R)$-algebra $(A, \pi)$ with the base
$(\mathcal R, {\mathfrak{M}})$. Let $(\mathcal R', {\mathfrak{M}}_{\mathcal R'})$ be an extension  of $(\mathcal R, {\mathfrak{M}})$ by an
$\mathcal R$-module $M$, where $\mbox{dim}_{\mathcal R}M<\infty$. In this section
we consider the problem of extending the given deformation $\lambda $
to a deformation with base $(\mathcal R', {\mathfrak{M}}_{\mathcal R'}).$

First let us consider the case of $1$-dimensional extension. Let
$$0\longrightarrow \bK\stackrel{i}{\longrightarrow}
\mathcal R' \stackrel{p}{\longrightarrow} \mathcal R \longrightarrow 0$$ be any
$1$-dimensional extension of $\mathcal R$. By the above proposition the isomorphism classes of $1$-dimensional extensions
of $\mathcal R$ are in one-to-one correspondence with the Harrison cohomology
$H_{Harr}^2(\mathcal R; \bK)$ of $\mathcal R$ with coefficients in $\bK$ where the $\mathcal R$ module structure on $\bK$ is given by $rk =\epsilon (r)k.$

Let $[f] \in H_{Harr}^2(\mathcal R; \bK)$. Suppose
$0\longrightarrow {\bK}\stackrel{i}{\longrightarrow} \mathcal R'
\stackrel{p}{\longrightarrow} \mathcal R \longrightarrow 0$
is a representative of the class of $1$-dimensional extensions of $\mathcal R$,
corresponding to the cohomology class $[f]$.

Let us recall how the algebra structure on $\mathcal R'$ is related to $f$.
Fix a splitting $q: \mathcal R \longrightarrow \mathcal R'$.
Let $\widehat{\epsilon}= \epsilon \circ p: \mathcal R'\longrightarrow \bK$ denote the
augmentation of $\mathcal R'$.
Then $b\mapsto (p(b), i^{-1}(b-q\circ p(b)))$ is a $\bK$-module
isomorphism $\mathcal R'\cong \mathcal R\oplus \bK.$ Let $(r,k)_q \in \mathcal R'$ denote the
inverse of $(r,k) \in \mathcal R\oplus \bK$ under the above isomorphism. The
cocycle $f$ representing the extension is determined by
$f(r_1, r_2) = i^{-1}((r_1,0)_q(r_2,0)_q -(r_1r_2,0)_q).$ On the
other hand, $f$ determines the algebra structure on $\mathcal R'$ by
$$(r_1,k_1)_q \cdot (r_2,k_2)_q := (r_1r_2, r_1\cdot k_2+ r_2\cdot k_1
+f(r_1,r_2))_q.$$
As in Section\,3, let $\{m_i\}_{i=1}^r$ be a fixed basis of the maximal ideal
$\mathfrak{M}$ of $\mathcal R$ with the dual basis $\{\xi _i\}_{i=1}^r$. Let $\psi _i = \alpha_{\lambda, \xi_i}, 1\leq i \leq r$ be the $1$-cochain introduced as in Section\,3.
Then by Proposition \ref{MC}, the deformation $\lambda $
can be written as
$$\lambda =
\pi + \sum_{i=1}^r m_i \otimes \psi_i ;
$$ where $\sum_{i=1}^r m_i \otimes \psi_i$ is a Maurer-Cartan element in $\mathfrak{M}\otimes \mathfrak{g}^\pi$.

Let $\{n_j\}_{1\leq j\leq r+1}$ be defined by $n_j = (m_j, 0)_q$ for $1\leq j \leq r$ and $n_{r+1} =
(0,1)_q.$ 
Then $\{n_j\}_{1\leq j\leq r+1}$ is a basis of the maximal ideal ${\mathfrak{M}}_{\mathcal R'}$
of $\mathcal R'$.

 A deformation $\Gamma$ with base $\mathcal R'$ extending $\lambda $ is entirely
determined by the following two facts:\\
$\bullet$  $\Gamma$ defined on $\g$ can be extended to the
category of $\mathcal R'$-modules, and \\
$\bullet$  if $\hat{\Gamma}$ is the unique extension of $\Gamma$ then $\hat{\Gamma}(r) = 0$ for every
$r\in R.$\\
Let $\psi \in \mathfrak{g}^1$ be any cochain. Define
\begin{equation}\label{gamma}
\Gamma(1_{\mathcal R'}\otimes \mu)
=  \pi(\mu) + \sum_{j=1}^r n_j \otimes \psi_j(\mu)+ n_{r+1}\otimes \psi(\mu),
\end{equation} where $\mu \in \g.$

Extending to the category of $\mathcal R'$-modules this defines a $\mathcal R'$-linear map

$$\Gamma: \mathcal{R'}\otimes \g \longrightarrow \mathcal{R'}\otimes \mbox{End}_A.$$

By universal property of free operad, we extend $\Gamma$ to a morphism
of operads $\widetilde{\Gamma}: {\mathcal F}(E_{\mathcal R'}) \longrightarrow \End(A_{\mathcal R'})$. Now, $\Gamma $ induces a $\cP_{\mathcal R'}$-algebra structure on
$A_{\mathcal R'}$ if and only if $\widetilde{\Gamma}((R))=0$ and it is clear from our
construction of $\Gamma $ that it extends the given deformation
$\lambda$.

As $\lambda$ defines an algebra structure on $\mathcal R \otimes A$, $[\lambda, \lambda]=0$. This implies (Proposition\,12.2.6 of \cite{LV})
\begin{equation}\label{lambda2}
\sum_{i=1}^r m_i \otimes \partial_{\pi}(\psi_i) + 1/2 \sum_{i,j=1}^r m_i m_j\otimes [\psi_i, \psi_j]=0.\end{equation}
Note that  $\Gamma$ induces an algebra structure on $\mathcal{R}'\otimes A$ is equivalent to saying that as an element of the Lie algebra $\mathfrak{g}_{\mathcal{R}'\otimes \cP, \mathcal{R}'\otimes A}
\cong\mathcal{R}' \otimes \mathfrak{g}_{\cP,A}$, $[\Gamma, \Gamma]$ vanishes. From the expression of $\Gamma$ we get
$$\begin{array}{ll}
&[\Gamma, \Gamma]\\
= &[\pi+ \sum_{i=1}^r n_i \otimes\psi_i + n_{r+1}\otimes \psi, \pi+ \sum_{i=1}^r n_i\otimes\psi_i + n_{r+1}\otimes\psi]\\
=& [\pi,\pi]+ 2\sum_{i=1}^r n_i \otimes \partial_{\pi} \psi_i + 2n_{r+1}\otimes \partial_{\pi} \psi + \sum_{i,j=1}^r n_i n_j \otimes [\psi_i, \psi_j]\\
&+2\sum_{i=1}^r n_i n_{r+1}\otimes [\psi_i, \psi] + n_{r+1}^2\otimes [\psi, \psi]
\end{array}$$ using the fact that $\partial_{\pi}(f)= [\pi, f]$.

We note that $[\pi, \pi]=0$. As $n_i n_{r+1}=0$ for $1\leq i\leq r$ ($\mathfrak{M}\bK=0$) and $n_{r+1}^2=0$ the above expression is equal to
$$2\sum_{i=1}^r n_i \otimes \partial_{\pi}\psi_i + 2n_{r+1}\otimes \partial_{\pi} \psi+ \sum_{i,j=1}^r n_i n_j \otimes [\psi_i, \psi_j].$$ 

Thus $\Gamma$ defines a $\mathcal{R}'\otimes \cP$-algebra structure on $\mathcal{R}'\otimes A$ extending $\lambda$, the $\mathcal{R}\otimes \cP$-algebra structure on $\mathcal{R}\otimes A$, if and only if
$$\begin{array}{ll}
&[\Gamma, \Gamma]=0\\
\Leftrightarrow & n_{r+1}\otimes \partial_{\pi} \psi= - \sum_{i=1}^r n_i \otimes \partial_{\pi} \psi_i- 1/2\sum_{i,j=1}^r n_i n_j \otimes [\psi_i, \psi_j]\\
\Leftrightarrow & (0, 1)_q\otimes \partial_{\pi}(\psi) \\
=& -\sum_{i=1}^r (m_i, 0)_q\otimes \partial_{\pi}(\psi_i) -1/2 \sum_{i, j=1}^r (m_i m_j, f(m_i, m_j))_q\otimes [\psi_i, \psi_j]\\& ~~(\mbox{using~ the~ isomorphism~ between~} \mathcal R' \cong \mathcal R\oplus \bK )\\
\Leftrightarrow &\partial_{\pi}(\psi)= -1/2\sum_{i, j=1}^r f(m_i, m_j)[\psi_i, \psi_j].
\end{array}$$

Let us define a $2$-cochain $\Phi$ on $\mathfrak{g}_{\cP, A}$ as follows: \\
$$\Phi= \sum_{i, j=1}^r f(m_i, m_j) [\psi_i, \psi_j].$$ This cochain is called the {\it obstruction cochain}.

\begin{prop}
The obstruction cochain $\Phi$ is a $2$-cocycle in $\mathfrak{g}_{\cP, A}$.
\end{prop}
{\bf Proof.} 

From equation \ref{lambda2}, 
\begin{equation}\label{m_i}
\sum_{i=1}^r(m_i, 0)_q \otimes \partial_\pi \psi_i +1/2 \sum_{i,j=1}^r (m_i m_j, 0)_q\otimes [\psi_i, \psi_j]=0
\end{equation}
as an element of $\mathcal R' \otimes \mathfrak{g}^\pi$.

Let $\partial_\pi'$ is the differential of $\mathcal R'\otimes \mathfrak{g}$ induced by $\partial_\pi$ in $\mathfrak{g}^\pi$. Then
$$\begin{array}{ll}
&\partial_\pi '(\sum_{i,j=1}^r n_in_j \otimes [\psi_i, \psi_j])\\
=&\sum_{i,j=1}^r n_i n_j\otimes \partial_\pi[\psi_i, \psi_j]\\
=&\big(\sum_{i,j=1}^r n_in_j\otimes \big([\partial_\pi \psi_i, \psi_j]- [\psi_i, \partial_\pi \psi_j]\big)\big)\\
=&\big([\sum_{i=1}^r n_i \otimes \partial_\pi \psi_i, \sum_{j=1}^r n_j\otimes \psi_j]- [\sum_{i=1}^r n_i \otimes \psi_i, \sum_{j=1}^r n_j\otimes \partial_\pi\psi_j]\big)\\
=&[\sum_{i=1}^r(m_i, 0)_q \otimes \partial_\pi \psi_i, \sum_{j=1}^r(m_j, 0)_q \otimes \psi_j] \\
& ~~~-[\sum_{i=1}^r(m_i, 0)_q \otimes  \psi_i, \sum_{j=1}^r(m_j, 0)_q \otimes \partial_\pi\psi_j]\\
=&-1/2 [\sum_{k, l=1}^r (m_km_l, 0)_q \otimes [\psi_k, \psi_l], ~\sum_{j=1}^r(m_j, 0)_q \otimes \psi_j]\\
&~~~+1/2 [\sum_{i=1}^r(m_i, 0)_q \otimes  \psi_i, ~\sum_{k, l=1}^r (m_km_l, 0)_q \otimes [\psi_k, \psi_l]]\\
&~~~~~~\mbox{using equation(\ref{m_i})}\\
=&-1/2\sum_{i,j,k=1}^r (m_i m_j m_k, f(m_im_j, m_k))_q\otimes [[\psi_i, \psi_j], \psi_k]\\
&~~~+ 1/2\sum_{i,j,k=1}^r( m_im_jm_k, f(m_i, m_j m_k))_q\otimes [\psi_i, [\psi_j, \psi_k]]\\
=&1/2\sum_{i, j, k=1}^r n_i n_j n_k\otimes ([\psi_i, [\psi_j, \psi_k]]- [[\psi_i, \psi_j], \psi_k]) \\
&~~~~~\mbox{as}~ f(m_i m_j, m_k)= f(m_i, m_jm_k), f~ \mbox{being a}~ 2\mbox{-cocycle}~\mbox{and}~ \mathfrak{M}\bK=0\\
=&1/2\sum_{i,j,k=1}^r n_in_jn_k \otimes [[\psi_i, \psi_k], \psi_j]\\
=& 0, ~~\mbox{by Jacobi identity}.
\end{array}
$$
On the other hand, $$
\partial_\pi '(\sum_{i,j=1}^r n_in_j \otimes [\psi_i, \psi_j])
  =\sum_{i,j=1}^r (m_im_j, f(m_i, m_j))_q \otimes \partial_\pi[\psi_i, \psi_j].$$
Hence, by previous argument, $$\sum_{i,j=1}^r (m_im_j, f(m_i, m_j))_q \otimes \partial_\pi[\psi_i, \psi_j]=0.$$
By equation\,(\ref{lambda2}), $\sum_{i,j=1}^r m_im_j\otimes \partial_\pi [\psi_i, \psi_j]=0$.
Hence, $$\partial_\pi\Phi= \sum_{i,j=1}^r f(m_i, m_j)\partial_\pi[\psi_i, \psi_j]=0.$$ \qed

The above proposition enables us to define a map from the set of $2$-cocycles $Z^2(\mathcal R, \bK)$ to $ H_\cP^2(A)$. We shall see that this map passes on to the cohomology.

Let $f$ and $f_1$ be two $2$-cocycles determining the same cohomology class in $ H_{Harr}^2(\mathcal R;\bK)$, that is $f-f_1= \delta h$ for some $1$-cochain $h \in Ch^1(\mathcal{R}, \bK)$. The obstruction cocycle for extending $\lambda$ to an extension corresponding to $f$ is given by $\sum_{i,j} f(m_i, m_j)[\psi_i, \psi_j]$. The obstruction cocycle corresponding to the cocycle $f_1$ is given by $\sum_{i,j} f_1(m_i, m_j)[\psi_i, \psi_j]$. Now,
$$\begin{array}{ll}
&\sum_{i, j} f(m_i, m_j)[\psi_i, \psi_j]- \sum_{i,j} f_1(m_i, m_j)[\psi_i, \psi_j]\\
=&\sum_{i, j} \delta h(m_i, m_j)[\psi_i, \psi_j]\\
=&\sum_{i, j}(m_i h(m_j) - h(m_im_j) + h(m_i)m_j) [\psi_i, \psi_j]\\
=&-\sum_{i, j} h(m_i, m_j)[\psi_i, \psi_j].
\end{array}
$$ The terms  $m_ih(m_j)$ and $h(m_i)m_j$ vanish as $\epsilon(\mathfrak{M})=0$.
On the other hand, from equation (\ref{lambda2}), 
$$\sum_{i=1}^r m_i \otimes \partial_{\pi}(\psi_i) =-1/2\sum_{i,j=1}^r m_i m_j\otimes [\psi_i, \psi_j].$$
Hence, $$
\begin{array}{ll}
&(h\otimes \mbox{id})\sum_{i=1}^r m_i \otimes \partial_\pi \psi_i= -1/2 \sum_{i,j=1}^r h(m_i m_j)[\psi_i, \psi_j]\\
\Rightarrow&\partial_{\pi}\big(2\sum_{i=1}^r h(m_i) \psi_i\big)= \sum_{i=1}^r h(m_i)\partial_{\pi}(\psi_i)= -\sum_{i, j=1}^r h(m_i m_j)[\psi_i, \psi_j].
\end{array}$$

The above consideration defines a map
$$\theta_{\lambda}:  H_{Harr}^2(\mathcal R; \bK) \longrightarrow  H_\cP^2(A),~\mbox{by}~
\theta_{\lambda}([f])=[\Phi],$$ where $[\Phi]$ is the cohomology class of $\Phi$. The map $\theta_{\lambda}$ is called the {\it obstruction map}. Proof of the following proposition is straightforward.
\begin{prop}
Let $\lambda$ be a deformation of a $\cP$-algebra $A$ with base $\mathcal R$ and let $\mathcal R'$ be a $1$-dimensional extension of $\mathcal R$ corresponding to the cohomology class $[f] \in H^{2}_{Harr}(\mathcal R;\bK)$. Then $\lambda$ can be extended to a deformation of $A$ with base $\mathcal R'$ if and only if the obstruction $\theta_{\lambda}([f])=0$.\qed
\end {prop}

We state the following proposition, proof of which is similar to the proof of corollary \,5.8 in \cite{fmm}.
\begin{prop}\label{d is onto}
Suppose that for a deformation $\lambda$ of a $\cP$-algebra $A$ with base $\mathcal R$, the differential $d\lambda: T\mathcal R \longrightarrow \mathbb{H}$ is onto.
Then the group of automorphisms $\mathcal{A}$ of the extension \begin{equation}\label{1-dimensional ext}
0\longrightarrow
\bK\stackrel{i}{\longrightarrow} \mathcal R'
\stackrel{p}{\longrightarrow}\mathcal R\longrightarrow 0
\end{equation} operates transitively on the set of equivalence classes of deformations $\Gamma$ of $A$ with base $\mathcal R'$ such that $p_{*}\Gamma=\lambda$. In other words, if $\Gamma$ exists, it is unique up to an isomorphism and an automorphism of this extension. \qed
\end{prop}

More generally, for an extension
$$0\longrightarrow
{M}\stackrel{i}{\longrightarrow} \mathcal R'
\stackrel{p}{\longrightarrow}\mathcal R\longrightarrow 0$$ with $\mbox{dim}_{\mathcal R}M< \infty$, the above arguments can be generalised.

The obstruction map in this more general situation is defined by
$$\begin{array}{cccc}
 \theta_{\lambda}:& H_{Harr}^2(\mathcal R; M) &\longrightarrow & M \otimes H_\cP^2(A)\\
& [f]& \mapsto & [\Phi] .
\end{array}$$
Then, as in the case of $1$-dimensional extension, we have the following.

\begin{prop}\label{for any extension}
 Let $\lambda$ be a deformation of a $\cP$-algebra $A$ with base $(\mathcal R ,\mathfrak{M})$ and let $M$ be a finite dimensional $\mathcal R$-module with $\mathfrak{M}M=0$. Consider an extension $\mathcal R'$ of $\mathcal R$
 $$0\longrightarrow
{M}\stackrel{i}{\longrightarrow} {\mathcal R'}
\stackrel{p}{\longrightarrow}\mathcal R\longrightarrow 0$$
corresponding to some $[f] \in H^{2}_{Harr}(\mathcal R;M)$. A deformation $\Gamma$ of $A$ with base $\mathcal R'$ such that $p_{*}\Gamma=\lambda$ exists if and only if the obstruction $\theta_{\lambda}([f])=0$.
If $d\lambda: T\mathcal R \longrightarrow \mathbb{H}$ is onto, then the deformation $\Gamma$, if it exists, is unique up to an isomorphism and an automorphism of the above extension. \qed
\end{prop}
We end this section with the following naturality property of the obstruction map, proof of which is similar to Proposition\,5.10 in \cite{fmm}.

\begin{prop}\label{obstructions are same}
Suppose  $\mathcal R_1$ and $\mathcal R_2$ are finite dimensional unital local algebras with augmentations $\varepsilon_1$ and $\varepsilon_2$, respectively.  Let $\phi : \mathcal R_{2} \longrightarrow \mathcal R_{1}$ be an algebra  homomorphism with $\phi(1)=1$ and  ${\varepsilon}_1 \circ \phi=\varepsilon_2$. Suppose $\lambda_2$ is a deformation of a $\cP$-algebra $A$ with base $\mathcal R_2$ and $\lambda_1 = \phi_{*} \lambda_2$ is the push-out via $\phi$. Then  the following diagram commutes.

$$\xymatrix{
&H_{Harr}^2({\mathcal R}_1; M) \ar[dd]^{\phi^*}\ar[dr]^{\theta_{\lambda_1}}&\\
&& H_{\cP}^2(A)\\
&H_{Harr}^2({\mathcal R}_2; M)\ar[ur]_{\theta_{\lambda_2}}&
}$$
$$Figure~1.$$\qed
\end {prop}
\medskip

\section {Construction of a Versal Deformation}

We begin the last section with the definition of the notion of a versal deformation of a $\cP$-algebra. The importance of versal deformation lies in the fact that it includes information of all other non-equivalent deformations of a given object. We use the results developed in the last two sections to give a constructive proof of existence of versal deformations of a given $\cP$-algebra. In fact, starting with the infinitesimal deformation as introduced in Section\,4,  we give an explicit construction of a versal deformation of a $\cP$-algebra by an inductive argument, using obstruction theory developed in Section\,5.

\begin{definition}
A formal deformation $\eta$ of a $\cP$-algebra $A$ with base $\mathcal R'$ is
called {\it versal}, if it satisfies the conditions below:\\
(i) for any formal deformation $\lambda$ of $A$ with base $\mathcal R$ there exists a homomorphism
$f : {\mathcal R'}\longrightarrow {\mathcal R}$ such that the deformation $\lambda$ is equivalent to $f_*\eta$;\\
(ii) if $\mathcal R$ satisfies the condition ${\mathfrak{M}}^2 = 0$, then $f$ is unique.
\end{definition}

We proceed to construct a versal deformation of a $\cP$-algebra $A$. As before, let $\mathbb{H}= H_\cP^1(A)$ and assume that $dim(\mathbb{H})<\infty$.  Let $\eta_1$ be the universal infinitesimal deformation with base $\mathcal C_1$ as constructed in Section \ref{universal}. Suppose for some $k\geq 1$ we have constructed a finite dimensional local algebra $\mathcal C_k$ and a  deformation $\eta_k$ of $A$ with base $\mathcal C_k$.
Let
$$\mu:H_{Harr}^2(\mathcal C_k;\bK)\longrightarrow (\mathcal Ch_{2} (\mathcal C_k))^\prime$$
be a homomorphism sending a cohomology class to a cocycle representing the class. Let
$$f_{\mathcal C_k}:Ch_{2} (\mathcal C_k) \longrightarrow H_{Harr}^2(\mathcal C_k;\bK)^\prime$$
 be the dual of $\mu$. By Proposition \ref{cohomology class corresponds to extension} ({\it ii}) we have the following extension of $\mathcal C_k$:
\begin{equation}\label{universal extension}
0\longrightarrow
H_{Harr}^2(\mathcal C_k;\bK)^\prime \stackrel{\bar{i}_{k+1}}{\longrightarrow}\bar{\mathcal C}_{k+1}\stackrel{\bar{p}_{k+1}}{\longrightarrow} {\mathcal C}_k \longrightarrow 0.
\end{equation}
The corresponding obstruction $\theta_{\eta_k}([f_{\mathcal C_k}]) \in H_{Harr}^2(\mathcal C_k;\bK)^\prime \otimes H_{\cP}^2(A)$ gives a linear map
$\omega_k:H_{Harr}^2(\mathcal C_k;\bK) \longrightarrow H_\cP^2(A)$
with the dual map
$${\omega_k}^\prime:H_\cP^2(A)^\prime \longrightarrow H_{Harr}^2(\mathcal C_k;\bK)^\prime .$$
We have an induced extension
$$ 0\longrightarrow coker (\omega'_{k})\longrightarrow \bar{\mathcal C}_{k+1}/\big(\bar{i}_{k+1}\circ \omega'_{k}(H_\cP^2(A)')\big) \longrightarrow \mathcal C_k \longrightarrow 0.$$
Since $coker (\omega'_k)\cong (ker (\omega_k))^\prime$,
it yields an extension
\begin{equation}\label{yields an extension}
0\longrightarrow (ker(\omega_k))^\prime \stackrel{i_{k+1}}{\longrightarrow} \mathcal C_{k+1}
\stackrel{p_{k+1}}{\longrightarrow} \mathcal C_k \longrightarrow 0
\end{equation}
where $\mathcal C_{k+1}= \bar{\mathcal C}_{k+1}/\big(\bar{i}_{k+1}\circ~ \omega_{k}^\prime (H_\cP^2(A)')\big)$ and  $i_{k+1}$, $p_{k+1}$ are the mappings induced by $\bar{i}_{k+1}$ and $\bar{p}_{k+1}$, respectively.
Observe that the algebra $\mathcal C_{k+1}$ is also local.  Since $\mathcal C_k$ is finite dimensional, the cohomology group $H_{Harr}^2(\mathcal C_k;\bK)$ is also finite dimensional and hence $\mathcal C_{k+1}$ is finite dimensional as well.

\begin{remark}\label{specific extension}
It follows from Proposition \ref{coefficients in M} that the specific extension (\ref{universal extension}) has the following ``universality property''. For any $\mathcal C_k$-module $M$ with $\mathfrak{M}M=0$, (\ref{universal extension}) admits a unique morphism into an arbitrary extension of $\mathcal C_k$:
$$0\longrightarrow M \longrightarrow \mathcal R' \longrightarrow \mathcal C_k \longrightarrow 0.$$
\end{remark}

\begin{prop}
The deformation $\eta_k$ with base $\mathcal C_{k}$ of a $\cP$-algebra $A$ admits an extension to a deformation with base $\mathcal C_{k+1}$, which is unique up to an isomorphism and an automorphism of the extension
$$0\longrightarrow
 (ker(\omega_k))^\prime \stackrel{i_{k+1}}{\longrightarrow} \mathcal C_{k+1}
\stackrel{p_{k+1}}{\longrightarrow} \mathcal C_k \longrightarrow 0.$$
\end{prop}
{\bf Proof.} From the above construction of the extension in \ref{yields an extension}, it is clear that $\theta_{\eta_k}([f_{{\mathcal C}_k}])=\omega_k|_{ker(\omega_k)}=0$. Therefore the proof is complete by proposition \ref{for any extension}.\qed

By induction the above process 
gives a sequence of finite dimensional local algebras $C_k$ and the deformation $\eta_k$ of the $\cP$-algebra $A$ with base $C_k$:
$$\bK\stackrel{p_1}{\longleftarrow}{\mathcal C}_1\stackrel{p_2}{\longleftarrow}\cdots \stackrel{p_k}{\longleftarrow}{\mathcal C}_k\stackrel{p_{k+1}}{\longleftarrow}{\mathcal C}_{k+1}$$ such that $(p_{k+1})_*\eta_{k+1}= \eta_k$.

As a consequence, we obtain a formal deformation $\eta$ of the $\cP$-algebra $A$ with base $\mathcal C= \ilim_{k\rightarrow \infty} \mathcal C_k$.

Proof of the following theorem is analogous to the proof of theorem\, 6.8 of \cite{fmm}.

\begin{thm}
Let $A$ be a $\cP$-algebra with $dim(\mathbb{H})<\infty$. Then the formal deformation $\eta$ with base $\mathcal C$ constructed above is a versal deformation of $A$.
\end{thm}

{\bf Acknowledgement} We express our sincere gratitude to the anonymous referee for his valuable comments which helped us to improve an earlier version of this paper.


\begin{thebibliography}{99}
\medskip


\bibitem{bal2} D.\ Balavoine, Deformations of algebras over a quadratic operad, \emph{Contemp. Math.} \textbf{202} (1997), 207-234.

\bibitem{bjt} H.-J.\ Baues, M.\ Jibladze and A.\ Tonks, Cohomology of monoids in monoidal categories, \emph{Contemp. Math.}, \textbf{202} (1997), 137-165.

\bibitem{fi1}A.\ Fialowski, Deformations of Lie algebras, \emph{Math.
USSR-Sbornyik} \textbf{55} (1986), pp.467-473.

\bibitem{ff}A.\  Fialowski and D.B. Fuchs, Construction of versal
deformation of Lie algebras, \emph{Journal of Funct. Anal.} \textbf{161}
(1999), pp. 76-110.

\bibitem{fp}A.\ Fialowski and M. Penkava, Deformation Theory of Infinity
Algebras, \emph{J. Algebra} \textbf{255} (2002), pp. 59-88.

\bibitem{fmm}A.\ Fialowski, A. Mandal and G. Mukherjee, Versal
deformations of Leibniz Algebras, \emph{J. K-Theory} \textbf{3} (2009), no. 2, 327-358.

\bibitem{ger}M.\ Gerstenhaber, The cohomology structure of an associative ring, \emph{Ann. Math.} \textbf{78} (1963), 267-288.

\bibitem{ger2}M.\ Gerstenhaber, On the deformation of rings and algebras, \emph{Ann. Math.} \textbf{79} (1964), 59-103.

\bibitem{gj} E. \ Getzler and J. D. S \ Jones, Operads, homotopy algebras and iterated integrals for double loop spaces, hep-th/9403055 (1994).

\bibitem{gk}V.\ Ginzburg and M.\ M.\ Kapranov, Koszul duality for operads, emph{Duke Math. J.} \textbf{76} (1994), 203-272.

\bibitem{harrison}D.\ K.\ Harrison, Commutative algebras and cohomology, \emph{Trans.\ Amer.\ Math.\ Soc.} \textbf{104} (1962), 191-204.


\bibitem{ks}M.\ Kontsevich and Y.\ Soibelman, Deformations of algebras over operads and the Deligne conjecture, Conf\'{e}rence Mosh\'{e} Flato 1999, Vol. II (Dijon), pp.\ 255-307, Math.\ Phys.\ Stud.\ \textbf{22}, Kluwer, 2000.

\bibitem{LV}J.-L.\ Loday and B. Vallette, in: \emph{Algebraic operads},
Grundlehren\ der\ mathematischen\ Wissenschaften, 346 Springer,
2012.

\bibitem{markl}M.\ Markl, Models for operads, \emph{Comm.\ Alg.}\ \textbf{24} (1996), 1471-1500.

\bibitem{mss}M.\ Markl, S.\ Shnider, and J.\ Stasheff, \emph{Operads in Algebra, Topology and Physics, Math.\ Surveys and Monographs} \textbf{96}, Amer.\ Math.\ Soc., 2002.

\bibitem{may1}J.\ P.\ May, The geometry of iterated loop spaces, \emph{Lectures Notes in Math.}\ \textbf{271}, Springer-Verlag, 1972.

\bibitem{may2}J.\ P.\ May, Definitions: operads, algebras and modules, in:  \emph{Operads: Proceedings of Renaissance Conferences} (Hartford, CT/Luminy, 1995), pp.\ 1-7, Contemp.\ Math. \textbf{202}, Amer.\ Math.\ Soc., Providence, RI, 1997.

\bibitem{sch}M. Schlessinger, Functors of Artin Rings,
\emph{Trans. Amer. Math. Soc.}\ 130 (1968), pp. 208-222.


\end{thebibliography}
\end{document}